\documentclass{amsart}
\usepackage{latexsym,amsmath,amssymb}

\theoremstyle{definition}
\theoremstyle{remark}

\numberwithin{equation}{section}

\begin{document}

\title[CONDITIONAL TYPE OPERATORS]
{MULTIPLICATION CONDITIONAL EXPECTATION TYPE OPERATORS ON ORLICZ  SPACES}

\author{\sc\bf Y. Estaremi }
\address{\sc Department of Mathematics, Payame Noor University, p. o. box: 19395-3697, Tehran, Iran.}
\email{estaremi@gmail.com, yestaremi@pnu.ac.ir}

\address{}

\thanks{}

\thanks{}

\subjclass[]{47B47}
\begin{abstract}
In this paper we consider a generalized conditional-type
Holder-inequality and investigate some classic properties of
multiplication conditional expectation type operators on
Orlicz-spaces.
\end{abstract}

\keywords{Conditional expectation- Orlicz space-
Holder-inequality- Compact operator- Essential norm.\\
This work has been done under supervision of Professor Ben De
Pagter, when the author has been in Delft University of Technology
for a six month visit
 }

\date{}

\dedicatory{}

\commby{}


\maketitle

\section{ \sc\bf Introduction }

\vspace*{0.3cm} Let $(\Omega, \Sigma, \mu)$ be a measure space and
$\mathcal{A}\subseteq \Sigma$ be sub $\sigma-$algebra. For a
sub-$\sigma$-finite algebra $\mathcal{A}\subseteq\Sigma$, the
conditional expectation operator associated with $\mathcal{A}$ is
the mapping $f\rightarrow E^{\mathcal{A}}f$, defined for all
non-negative $f$ as well as for all $f\in L^1(\Sigma)$ and $f\in
L^{\infty}(\Sigma)$, where $E^{\mathcal{A}}f$, by the
Radon-Nikodym theorem, is the unique $\mathcal{A}$-measurable
function satisfying
$$\int_{A}fd\mu=\int_{A}E^{\mathcal{A}}fd\mu, \ \ \ \forall A\in \mathcal{A} .$$
As an operator on $L^{1}({\Sigma})$ and $L^{\infty}(\Sigma)$,
$E^{\mathcal{A}}$ is idempotent and
$E^{\mathcal{A}}(L^{\infty}(\Sigma))=L^{\infty}(\mathcal{A})$ and
$E^{\mathcal{A}}(L^1(\Sigma))=L^1(\mathcal{A})$. Thus it can be
defined on all interpolation spaces of $L^1$ and $L^{\infty}$ such
as, Orlicz spaces\cite{besh}. If there is no possibility of
confusion, we write $E(f)$ in place of $E^{\mathcal{A}}(f)$. This
operator will play a major role in our work and we list here some
of its useful properties:

\vspace*{0.2cm} \noindent $\bullet$ \  If $g$ is
$\mathcal{A}$-measurable, then $E(fg)=E(f)g$.

\noindent $\bullet$ \ $\varphi(E(f))\leq E(\varphi(f))$, where
$\varphi$ is a convex function.

\noindent $\bullet$ \ If $f\geq 0$, then $E(f)\geq 0$; if $f>0$,
then $E(f)>0$.

\noindent $\bullet$ \ For each $f\geq 0$, $\sigma(f)\subseteq
\sigma(E(f))$.

\vspace*{0.2cm}\noindent A detailed discussion and verification of
most of these properties may be found in \cite{rao}. We recall
that an $\mathcal{A}$-atom of the measure $\mu$ is an element
$A\in\mathcal{A}$ with $\mu(A)>0$ such that for each
$F\in\mathcal{A}$, if $F\subseteq A$, then either $\mu(F)=0$ or
$\mu(F)=\mu(A)$. A measure space $(X,\Sigma,\mu)$ with no atoms is
called non-atomic measure space.\cite{z} \ \ \

Let $(\Omega, \Sigma, \mu)$ be a measure space and
$\mathcal{A}\subseteq \Sigma$ be sub $\sigma-$algebra, such that
$(\Omega, \mathcal{A}, \mu_{\mathcal{A}})$ has finite subset
property. $E^{\mathcal{A}}=E$ is conditional expectation with
respect to $\mathcal{A}$. It is well-known fact that every
$\sigma$-finite measure space $(\Omega, \Sigma,\mu)$ can be
partitioned uniquely as $\Omega=\left
(\bigcup_{n\in\mathbb{N}}C_n\right )\cup B$, where
$\{C_n\}_{n\in\mathbb{N}}$ is a countable collection of pairwise
disjoint $\Sigma$-atoms and $B$, being disjoint from each $C_n$,
is non-atomic.\cite{z}\ \ \

Operators in function spaces defined by conditional expectations
were first studied, among others, by S - T.C. Moy \cite{mo}, Z.
Sidak \cite{s} and H.D. Brunk \cite{b} in the setting of $L^p$
spaces. Conditional expectation operators on various function
spaces exhibit a number of remarkable properties related to the
underlying structure of the given function space or to the metric
structure when the function space is equipped with a norm. P.G.
Dodds, C.B. Huijsmans and B. de Pagter \cite{dhd} linked these
operators to averaging operators defined on abstract spaces
earlier by J.L. Kelley \cite{k}, while A. Lambert \cite{l} studied
their link to classes of multiplication operators which form
Hilbert $C^*$-modules. J.J. Grobler and B. de Pagter \cite{gd}
showed that the classes of partial integral operators, studied by
A.S. Kalitvin and others \cite{afkz, akn, akz, dou, kz}, were a
special case of conditional expectation operators. Recently, J.
Herron studied operators $EM_u$ on $L^p$
spaces in \cite{ her}. \\

Also, in\cite{e, ej} we investigate some classic properties of
multiplication conditional expectation operators $M_wEM_u$ on
$L^p$ spaces. In the present paper we continue the investigation
of some classic properties of the operator $EM_u$ on Orlicz spaces
by considering Generalized conditional-type Holder inequality. \ \
\\
 Let us now introduce the definition of convexity for functions
of n variables and later some particular criteria for convex
functions of $2$ variables.

\vspace*{0.3cm} {\bf Definition}\\

 Let be $f:\mathbb{R}^n \rightarrow \mathbb{R}$ and $x=(x_1, . . . , x_n)\in \mathbb{R}^n$,
we say that the function $f(x)$ is convex in $\mathbb{R}^n$ (or in
a subset of $\mathbb{R}^n$) if
$$f(\lambda x + (1 - \lambda)y)\leq \lambda f(x) + (1
- \lambda)f(y),$$ for any $x=(x_1, . . . , x_n)$ and $y=(x_1, . .
. , x_n)$ and for any $0 \leq \lambda \leq 1$.\\

It is well known fact that:

 Given a nice function $f$ of $2$ variables
in a set in the plane, the function $f$ is convex if and only if
both the following properties are true in such set:

(1) $ f_{x, x}(x, y)f_{y, y}(x, y)-(f_{x, y}(x, y))^2\geq0$,\\

(2) $f_{x, x}(x, y) \geq 0$ and $f_{y, y}(x, y) \geq 0$.\\

Let $\Phi:\mathbb{R}\rightarrow\mathbb{R}^{+}$ be a continuous
convex function such that\\

(1)$\Phi(x)=0$ if and only if $x=0$.\\

(2) $\Phi(x)=\Phi(-x)$.\\

(3) $\lim_{x\rightarrow\infty}\frac{\Phi(x)}{x}=\infty$, $\lim_{x\rightarrow\infty}\Phi(x)=\infty$.\\

 The function $\Phi$ is called Young's function. With each Young's
function $\Phi$, one can associate another
 convex function $\Psi:\mathbb{R}\rightarrow\mathbb{R}^{+}$ having similar properties, which is defined by
$$\Psi(y)=\sup\{x|y|-\Phi(x):x\geq0\}, \ \ y\in\mathbb{R}.$$
Then $\Psi$ is called complementary Young function to $\Phi$. A
Young function $\Phi$ is said to satisfy the $\bigtriangleup_{2}$
condition (globally) if $\Phi(2x)\leq k\Phi(x), \ x\geq x_{0}\geq0
(x_{0}=0)$ for some constant $k>0$. Also, $\Phi$ is said to
satisfy the $\bigtriangleup'(\bigtriangledown')$ condition,
 in symbols $\Phi\in \bigtriangleup'( \Phi\in \bigtriangledown')$, if $\exists c>0$  $(b>0)$ such that
$$\Phi(xy)\leq c\Phi(x)\Phi(y), \ \ \ x,y\geq x_{0}\geq 0$$(and
$$\Phi(bxy)\geq \Phi(x)\Phi(y), \ \ \ x,y\geq y_{0}\geq 0).$$
If $x_{0}=0(y_{0}=0)$, then these conditions are said to hold
globally. If $\Phi\in \bigtriangleup'$, then $\Phi\in
\bigtriangleup_{2}$.

 Let $\Phi_{1}, \Phi_{2}$ be two Young
functions, then $\Phi_1$ is stronger than $\Phi_2$,
$\Phi_1\succ\Phi_2$ [or $\Phi_2\prec\Phi_1$] in symbols, if
$$\Phi_2(x)\leq\Phi_1(ax), \ \ \ x\geq x_0\geq0$$
for some $a_0\geq0$ and $x_0$, if $x_0=0$ then this condition is
said to hold globally.

Let $\Phi$ is a Young function, then the set of
$\Sigma-$measurable functions
$$L^{\Phi}(\Sigma)=\{f:\Omega\rightarrow \mathbb{C}:\exists k>0,
\int_{\Omega}\Phi(kf)d\mu<\infty\}$$ is a Banach space, with
respect to the norm
$N_{\Phi}(f)=\inf\{k>0:\int_{\Omega}\Phi(\frac{f}{k})d\mu\leq1\}$.
$(L^{\Phi}(\Sigma), N_{\Phi}(.))$ is called Orlicz space.

Let $\Phi$ be a Young function and $f\in L^{\Phi}(\Sigma)$. Since $\Phi$ is convex, by Jensen's inequality $\Phi(E(f))\leq E(\Phi(f))$ and so

$$\int_{\Omega}\Phi(\frac{E(f)}{N_{\Phi}(f)})d\mu=\int_{\Omega}\Phi(E(\frac{f}{N_{\Phi}(f)}))d\mu
$$$$\leq\int_{\Omega}E(\Phi(\frac{f}{N_{\Phi}(f)}))d\mu=\int_{\Omega}\Phi(\frac{f}{N_{\Phi}(f)})d\mu\leq1.$$
This implies that $N_{\Phi}(E(f))\leq N_{\Phi}(f)$ i.e, $E$ is a contraction on Orlicz spaces.

\
\

We say that $(E, \Phi)$ satisfies in Generalized conditional- type Holder inequality, if there exist some positive constant $C$, such that for all $f\in L^{\Phi}(\Omega, \Sigma, \mu)$ and $g\in L^{\Psi}(\Omega, \Sigma, \mu)$ we have
$$E(|fg|)\leq C \Phi^{-1}(E(\Phi(|f|)))\Psi^{-1}(E(\Psi(|g|))),$$
where $\Psi$ is complementary Young function to $\Phi$.

\

In the sequel as Lemma 1.2 , Lemma 1.3 and Lemma 1.5 we give some
conditions for $E$ or $\Phi$ or $\Psi$ or jointly to get
Generalized conditional- type Holder inequality.

\vspace*{0.3cm} {\bf Lemma 1.1}

Let $\Phi$ and $\Psi$ be complementary Young functions. If there exist $C_{1}, C_{2}>0$ such that
$$E(\Phi(\frac{f}{\Phi^{-1}(E(\Phi(f)))}))\leq C_{1}, \ \ \ E(\Psi(\frac{g}{\Psi^{-1}(E(\Psi(g)))}))\leq C_{2}.$$

Then $(E, \Phi)$ satisfies in Generalized conditional-type Holder
inequality.

\vspace*{0.3cm} {\bf Proof} Let $f\in L^{\Phi}(\Omega, \Sigma, \mu)$ and $g\in L^{\Psi}(\Omega, \Sigma, \mu)$. By Young inequality we have

$fg\leq \Phi(f)+\Psi(g)$. If we replace $f$ to $\frac{f}{\Phi^{-1}(E(\Phi(f)))}$ and $g$ to $\frac{g}{\Psi^{-1}(E(\Psi(g)))}$. We have

$$\frac{fg}{\Phi^{-1}(E(\Phi(f)))\Psi^{-1}(E(\Psi(g)))}\leq \Phi(\frac{f}{\Phi^{-1}(E(\Phi(f)))})+\Psi(\frac{g}{\Psi^{-1}(E(\Psi(g)))}).$$
By taking $E$ we have

$$\frac{E(fg)}{\Phi^{-1}(E(\Phi(f)))\Psi^{-1}(E(\Psi(g)))}\leq E(\Phi(\frac{f}{\Phi^{-1}(E(\Phi(f)))}))+E(\Psi(\frac{g}{\Psi^{-1}(E(\Psi(g)))}))\leq C_{1}+C_{2}.$$
This implies that
$$E(|fg|)\leq C \Phi^{-1}(E(\Phi(|f|)))\Psi^{-1}(E(\Psi(|g|))),$$
where $C=C_{1}+C_{2}$.\\

By using some basic facts about closed convex sets in Banach
spaces in section I.2 and I.3 of \cite{et} we have the following
lemma.\\

\vspace*{0.3cm} {\bf Lemma 1.2}

A mapping $F:\mathbb{R}^n\rightarrow [0,\infty)$ can be written in
the form

$$F(x_1, x_2,...,x_n)=\inf_{a\in A}\Sigma^{n}_{i=1}a_ix_i$$

 for some
countable set $A\subseteq \mathbb{R}^n_+$ if and only if $F$ is
concave, lower semicontinuous and positive homogeneous.
\\

Also, we recall a generalization of Jensen-inequality that is
proved by Markus Haase in \cite{mh}.
\\

\vspace*{0.3cm} {\bf Theorem 1.3}

Let $(\Omega, \Sigma, \mu)$ and $(\Omega', \Sigma', \mu')$ be
measure spaces, let $C\subseteq \mathcal{M}(\Omega, \Sigma,
\mu)_+$ be a subcone and let $T:C\rightarrow \mathcal{M}(\Omega',
\Sigma', \mu')_+$ be a monotone, subadditive and positively
homogeneous operator. Let $F:\mathbb{R}^n_+\rightarrow [0,\infty)$
be given by
$$F(x_1, x_2,...,x_n)=\inf_{a\in A}\Sigma^{n}_{i=1}a_ix_i$$
for some countable set $A\subseteq \mathbb{R}^n_+$. Then if $f_1,
f_2,...,f_n\in C$ such that $F(f_1,..., f_n)$ is in $C$, one has

$$T[F(f_1,..., f_n)]\leq F(Tf_1,..., Tf_n)$$
as an inequality in $\mathcal{M}(\Omega', \Sigma', \mu')$.

\vspace*{0.3cm} {\bf Corollary 1.4}

Let $\Phi$ and $\Psi$ be complementary Young functions such that
$\Phi''(x)\Psi''(y)\Phi(x)\Psi(y)-(\Phi'(x)\Psi'(y))^2\geq0$ and
the function $F(x,y)=\Phi^{-1}(x)\Psi^{-1}(y)$ is positively
homogeneous on $\mathbb{R}^2_+$. Then $(E, \Phi)$ satisfies in
Generalized conditional-type Holder inequality, for every
conditional expectation operator $E$.

\vspace*{0.3cm} {\bf Proof} Since $\Phi$ and $\Psi$ are
continuous, then the map $F(x,y)=\Phi^{-1}(x)\Psi^{-1}(y)$ is also
continuous. By assumptions the map $F$ is concave, lower
semicontinuous and positive homogeneous. By replacing $T$ with
conditional expectation operator $E$ in Theorem 1.3 we have

$$E(\Phi^{-1}(f)\Psi^{-1}(g))\leq \Phi^{-1}(E(f))\Psi^{-1}(E(g))$$
for all nonnegative measurable function on measure space $(\Omega,
\Sigma, \mu)$. Direct computation shows that for all $f\in
L^{\Phi}(\Omega, \Sigma, \mu)$ and $g\in L^{\Psi}(\Omega, \Sigma,
\mu)$,

$$E(|fg|)\leq \Phi^{-1}(E(\Phi(|f|)))\Psi^{-1}(E(\Psi(|g|))).$$

\vspace*{0.3cm} {\bf Lemma 1.5}\\
 Let $E$ be the conditional
expectation operator and $\Phi$ and $\Psi$ be complementary Young
functions. If there exists positive constant $C$ such that
$E(fg)\leq C E(f)E(g)$, for positive measurable functions $f\in
L^{\Phi}(\Omega, \Sigma, \mu)$ and $g\in L^{\Psi}(\Omega, \Sigma,
\mu)$. Then $(E, \Phi)$
satisfies in Generalized conditional- type Holder inequality.\\

\vspace*{0.3cm} {\bf Proof} Let $f\in L^{\Phi}(\Omega, \Sigma,
\mu)$ and $g\in L^{\Psi}(\Omega, \Sigma, \mu)$ such that $f>0,
g>0$. Since $\Phi^{-1}$ and $\Psi^{-1}$ are concave, then
$$E(f)=E(\Phi^{-1}\Phi(f))\leq \Phi^{-1}(E(\Phi(f)), \ \ \ \ \ E(g)=E(\Psi^{-1}\Psi(g))\leq
\Psi^{-1}(E(\Psi(g)).$$

This implies that

$$E(fg)\leq CE(f)E(g)\leq
C\Phi^{-1}(E(\Phi(f))\Psi^{-1}(E(\Psi(g)).$$

So for all $f\in L^{\Phi}(\Omega, \Sigma, \mu)$ and $g\in
L^{\Psi}(\Omega, \Sigma, \mu)$ we have

$$E(|fg|)\leq C\Phi^{-1}(E(\Phi(|f|))\Psi^{-1}(E(\Psi(|g|)).$$

\vspace*{0.3cm} {\bf Example 1.6}\\
 (a) If
$\Phi(x)=\frac{x^p}{p}$, $1\leq p<\infty$. Then for all $f\in
L^{\Phi}(\Omega, \Sigma, \mu)=L^p(\Sigma)$ and $g\in
L^{\Psi}(\Omega, \Sigma, \mu)=L^q(\Sigma)$, where
$\frac{1}{p}+\frac{1}{q}=1$, we have

$$E(|fg|)\leq(E(|f|^p))^{\frac{1}{p}}(E(|g|^q))^{\frac{1}{q}}.$$

\
(b) Let $\Omega=[-1,1]$, $d\mu=\frac{1}{2}dx$ and $\mathcal{A}=<\{(-a,a):0\leq a\leq1\}>$ (Sigma algebra generated by symmetric intervals).
Then
 $$E^{\mathcal{A}}(f)(x)=\frac{f(x)+f(-x)}{2}, \ \ x\in \Omega,$$
 where $E^{\mathcal{A}}(f)$ is defined. Thus $E^{\mathcal{A}}(|f|)\geq\frac{|f|}{2}$. Hence $|f|\leq2E(|f|)$.
 Let $\Phi$ be a Young function. For each $f\in L^{\Phi}(\Omega, \Sigma, \mu)$ we have  $\Phi(|f|)\leq2E(\Phi(|f|))$. Since $\Phi^{-1}$ is also increasing and concave
 $$|f|=\Phi^{-1}(\Phi(|f|))\leq \Phi^{-1}(2E(\Phi(|f|)))\leq2\Phi^{-1}(E(\Phi(|f|))),$$
 so $|f|\leq2\Phi^{-1}(E(\Phi(|f|)))$. Similarly for $g\in L^{\Psi}(\Omega, \Sigma, \mu)$ we have
 $|g|\leq2\Psi^{-1}(E(\Psi(|g|)))$. Thus
 $$|fg|\leq 4 \Phi^{-1}(E(\Phi(|f|)))\Psi^{-1}(E(\Psi(|g|))).$$
 By taking $E$ we have
 $$E(|fg|)\leq 4 \Phi^{-1}(E(\Phi(|f|)))\Psi^{-1}(E(\Psi(|g|))).$$

 Also, by lemma 1.1 we have $C_{1}=\Phi(2)$, $C_{2}=\Psi(2)$ and
 $C=\Phi(2)+\Psi(2)$.\\

 (c) Let $dA(z)$ be the normalized Lebesgue measure on open unit disc $\mathbb{D}$.
 Recall that for $1\leq p<\infty$ the Bergman space $L^p_a(\mathbb{D})$
 is the collection of all functions $f\in H(\mathbb{D})$, holomorphic functions on
 $\mathbb{D}$, for which $\int_{\mathbb{D}}|f(z)|^pdA(z)<\infty$.
Let $\mathcal{A}$ be the $\sigma$-algebra generated by
$\{(z^n)^{-1}(U): U\subseteq\mathbb{C}\ \mbox{is open}\}$. Then
$$E(u)(\xi)=\frac{1}{n}\sum_{\zeta^n=\xi}u(\zeta),
\ \ u\in H(\mathbb{D}), \ \ \xi\in\mathbb{D}\setminus \{0\} ,$$
(see \cite{ho}).  Note that $|u|\leq nE(|u|)$. By the same method
of part(b), for every Young function $\Phi$, $(E, \Phi)$ satisfies
in Generalized conditional- type Holder inequality.\\

(d)Let $X=[0,1]$, $\Sigma$=sigma algebra of Lebesgue measurable
subset of $X$, $\mu$=Lebesgue measure on $X$. Fix $n\in\{2, 3, 4 .
. .\}$ and let $s:[0, 1]\rightarrow[0, 1]$ be defined by
$s(x)=x+\frac{1}{n}$(mod 1). Let $\mathcal{B}=\{E\in
\Sigma:s^{-1}(E)=E\}$. In this case
$$E^{\mathcal{B}}(f)(x)=\sum^{n-1}_{j=0}f(s^{j}(x)),$$ where $s^j$
denotes the jth iteration of $s$. The functions $f$ in the range
of $E^{\mathcal{B}}$ are those for which the n graphs of f
restricted to the intervals $[\frac{j-1}{n},\frac{j}{n}]$, $1\leq
j\leq n$, are all congruent. Also, $|f|\leq nE^{\mathcal{B}}(|f|)$
a.e. By the same method of part(b), for every Young function
$\Phi$, $(E, \Phi)$ satisfies in Generalized conditional-type
Holder inequality.\\

(e) Let $(\Omega, \Sigma, \mu)$ be a measure space and
$\mathcal{A}\subseteq \Sigma$ be sub $\sigma-$algebra. If there
exists $C_{0}>0$ such that $|f|\leq C_{0}E^{\mathcal{A}}(|f|)$.
Then $(E^{\mathcal{A}}, \Phi)$ satisfies in Generalized
conditional- type Holder inequality, for every Young function
$\Phi$.

The part (e) of the last example and proposition 2.2 of \cite{lam}
show that there are many conditional expectation operators $E$
such that the Generalized conditional-type Holder inequality holds
for $(E, \Phi)$.

\section{ \sc\bf BOUNDEDNESS AND COMPACTNESS OF $EM_{u}$ ON ORLICZ SPACES}

\vspace*{0.3cm} {\bf Theorem 2.1.} Let $T=EM_{u}:L^{\Phi}(\Omega, \Sigma, \mu)\rightarrow L^{0}(\Omega, \Sigma, \mu)$ such that
 $T(f)=E(uf)$ for $f\in L^{\Phi}(\Omega, \Sigma, \mu)$ is well defined, then the followings hold.

\vspace*{0.3cm}
(a) If $T$ is bounded on $L^{\Phi}(\Omega, \Sigma, \mu)$, then $E(u)\in L^{\infty}(\Omega, \mathcal{A}, \mu)$.

\vspace*{0.3cm}

(b) If $\Phi\in \bigtriangleup'$(globally) and $T$ is bounded on $L^{\Phi}(\Omega, \Sigma, \mu)$, then  $\Psi^{-1}(E(\Psi(u)))\in L^{\infty}(\mathcal{A})$.

\vspace*{0.3cm} (c) If $(E, \Phi)$ satisfies in Generalized
conditional-type Holder inequality and $\Psi^{-1}(E(\Psi(u)))\in
L^{\infty}(\mathcal{A})$, then $T$ is bounded.

\

In this case, $\|T\|\leq C\|\Psi^{-1}(E(\Psi(u)))\|_{\infty}$.

 \vspace*{0.3cm} {\bf Proof.} (a) Suppose that $E(u)\notin
L^{\infty}(\Omega, \mathcal{A}, \mu)$. If we set

$E_{n}=\{x\in \Omega:|E(u)(x)|>n\}$,  for all $n\in \mathbb{N}$,
 then $E_{n}\in \mathcal{A}$ and $\mu(E_{n})>0$.

 Since $(\Omega, \mathcal{A}, \mu)$  has finite subset property, we can assume that $0<\mu(E_{n})<\infty$, for all $n\in \mathbb{N}$.
 By definition of $E_{n}$ we have

 $$T(\chi_{E_{n}})=E(u\chi_{E_{n}})=E(u)\chi_{E_{n}}>n\chi_{E_{n}}.$$

 Since Orlicz's norm is monotone, then

 $$\|T(\chi_{E_{n}})\|_{\Phi}>\|n\chi_{E_{n}}\|_{\Phi}=n\|\chi_{E_{n}}\|_{\Phi}.$$

 This implies that $T$ isn't bounded.

(b)If $\Psi^{-1}(E(\Psi(u)))\notin L^{\infty}(\mathcal{A})$, then $\mu(E_{n})>0$. Where

  $$E_{n}=\{x\in \Omega:\Psi^{-1}(E(\Psi(u)))(x)>n\}$$

  and so $E_{n}\in \mathcal{A}$. Since $\Phi\in \bigtriangleup'$, then $\Psi\in\bigtriangledown'$, i.e., $\exists b>0$ such that

  $$\Psi(bxy)\geq\Psi(x)\Psi(y), \ \ \ x, y\geq0.$$
  Also, $\Phi\in \triangle_{2}$. Thus $(L^{\Phi})^{\ast}=L^{\Psi}$ and so $T^{\ast}=M_{\bar{u}}:L^{\Psi}(\mathcal{A})\rightarrow L^{\Psi}(\Sigma)$, is also bounded. Hence for each $k>0$ we have

  $$\int_{\Omega}\Psi(\frac{ku\chi_{E_{n}}}{N_{\Psi}(\chi_{E_{n}})})d\mu=\int_{\Omega}\Psi(ku\chi_{E_{n}}\Psi^{-1}(\frac{1}{\mu(E_{n})}))d\mu
  =\int_{\Omega}\Psi(ku\Psi^{-1}(\frac{1}{\mu(E_{n})}))\chi_{E_{n}}d\mu$$
  $$\geq\int_{E_{n}}\Psi(u)\Psi(\frac{k\Psi^{-1}(\frac{1}{\mu(E_{n})})}{b})d\mu
  \geq\left(\int_{E_{n}}E(\Psi(u))d\mu\right)\Psi(\frac{k}{b^2})\Psi(\Psi^{-1}(\frac{1}{\mu(E_{n})})$$

  $$\geq\Psi(n)\mu(E_{n})\frac{1}{\mu(E_{n})}\Psi(\frac{k}{b^2})=\Psi(n)\Psi(\frac{k}{b^2}).$$

  Thus
$$\int_{\Omega}\Psi(\frac{ku\chi_{E_{n}}}{N_{\Psi}(\chi_{E_{n}})})d\mu=\int_{\Omega}\Psi(kM_{u}(f_{n}))d\mu
\geq\Psi(n)\Psi(\frac{k}{b^2})\rightarrow\infty$$ as
$n\rightarrow\infty$, where
$f_{n}=\frac{\chi_{E_{n}}}{N_{\Psi}(\chi_{E_{n}})}$. Thus
$N_{\Psi}(M_{u}(f_{n}))\rightarrow\infty$, as
$n\rightarrow\infty$. This is a contradiction, since $M_{u}$ is
bounded.

(c) Put $M=\|\Psi^{-1}(E(\Psi(u)))\|_{\infty}$. For $f\in L^{\Phi}(\Omega, \Sigma, \mu)$ and $g\in L^{\Psi}(\Omega, \Sigma, \mu)$ we have

$$\int_{\Omega}\Phi(\frac{E(uf)}{CMN_{\Phi}(f)})d\mu=\int_{\Omega}\Phi(\frac{E(u\frac{f}{N_{\Phi}(f)})}{CM})d\mu$$

 $$\leq\int_{\Omega}\Phi(\frac{C \Phi^{-1}(E(\Phi(|\frac{f}{N_{\Phi}(f)})|)))\Psi^{-1}(E(\Psi(|u|)))}{CM})d\mu$$

 $$ \leq\int_{\Omega}\Phi(\Phi^{-1}(E(\Phi(|\frac{f}{N_{\Phi}(f)}|))))d\mu\leq\int_{\Omega}E(\Phi(\frac{f}{N_{\Phi}(f)}))d\mu$$

 $$=\int_{\Omega}\Phi(\frac{f}{N_{\Phi}(f)})d\mu\leq1.$$

 So $N_{\Phi}(E(uf))\leq CMN_{\Phi}(f)$. Thus $T=EM_{u}$ is bounded and $\|T\|\leq C\|\Psi^{-1}(E(\Psi(u)))\|_{\infty}$.

\vspace*{0.3cm} {\bf Corollary 2.2.} \\(a) If  $(E, \Phi)$
satisfies in Generalized conditional-type Holder inequality and
$\Phi\in \bigtriangleup'$(globally), then $T$ is bounded if and
only if
 $\Psi^{-1}(E(\Psi(u)))\in L^{\infty}(\mathcal{A})$.\\

 (b) If $\Psi\prec x$ and $(E, \Phi)$ satisfies in Generalized
conditional-type Holder inequality, then $T$ is bounded if and
only if $\Psi^{-1}(E(\Psi(u)))\in L^{\infty}(\mathcal{A})$.

\vspace*{0.3cm} {\bf Proof.} (b) Since $\Psi\prec x$
 then $EM_{\Psi(u)}\leq KEM_{u}$ for some $K>0$. If $\Psi^{-1}(E(\Psi(u)))\notin
 L^{\infty}(\mathcal{A})$, then the operator $EM_{\Psi(u)}$ is not
 bounded and so $T=EM_u$ is not bounded.

\vspace*{0.3cm} {\bf Theorem 2.3.} Let  $T=EM_{u}$ be bounded on
$L^{\Phi}(\Sigma)$, then the following hold. \ \

(a) If $T$ is compact, then
 $$N_{\varepsilon}(E(u))=\{x\in \Omega:E(u)(x)\geq\varepsilon\}$$
 consists of finitely many $\mathcal{A}-$atoms, for all  $\varepsilon>0$.
\
\

(b) If $T$ is compact and $\Phi\in \bigtriangleup'$(globally), then $N_{\varepsilon}(\Psi^{-1}(E(\Psi(u))))$

consists of finitely many $\mathcal{A}-$atoms, for all $\varepsilon>0$, where

$$N_{\varepsilon}(\Psi^{-1}(E(\Psi(u))))=\{x\in \Omega:\Psi^{-1}(E(\Psi(u))(x)\geq\varepsilon\}.$$
\
\

(c) If $(E, \Phi)$ satisfies in Generalized conditional-type
Holder inequality and $N_{\varepsilon}(\Psi^{-1}(E(\Psi(u)))$
consists of finitely many $\mathcal{A}-$atoms, for all
$\varepsilon>0$, then $T$ is compact.

\vspace*{0.3cm}

\vspace*{0.3cm} {\bf Proof.} (a) If there exists
$\varepsilon_{0}>0$, such that $N_{\varepsilon _0}(E(u))$ consists
of infinitely many $\mathcal{A}-$atoms or a non-atomic subset of
positive measure. Since $(\Omega, \mathcal{A}, \mu)$  has finite
subset property. In both cases, we can find a sequence of disjoint
$\mathcal{A}-$measurable subsets $\{A_{n}\}_{n\in \mathbb{N}}$ of
$N_{\varepsilon_{0}}(E(u))$ with $0<\mu(A_{n})<\infty$. Let
$f_{n}=\frac{\chi_{A_{n}}}{N_{\Phi}(\chi_{A_{n}})}$. Hence
 $$|f_{n}-f_{m}|=|f_{n}+f_{m}|=|f_{n}|+|f_{m}|$$
 for $n\neq m$. Also, $E(uf_{n})=E(u)f_{n}\geq \varepsilon_{0}f_{n}$. By monotonicity of $N_{\Phi}(.)$ we have

 $$N_{\Phi}(|E(uf_{n})-E(uf_{m})|)=N_{\Phi}(|E(u)(f_{n}-f_{m})|)$$
 $$=N_{\Phi}(|E(u)|(|f_{m}|+|f_{m}|))\geq N_{\Phi}(|E(u)|f_{n})\geq \varepsilon_{0}N_{\Phi}(f_{n})= \varepsilon_{0}.$$

 Thus $N_{\Phi}(|E(uf_{n})-E(uf_{m})|)\geq\varepsilon_{0}$. This implies that $T$ can not be compact.

\
\
\

(b) Suppose that, there exists $\varepsilon_{0}>0$, such that $N_{\varepsilon}(\Psi^{-1}(E(\Psi(u))))$ doesn't consist finitely many $\mathcal{A}-$atoms.
 Since $T$ is compact, so $T^{\ast}=M_{\bar{u}}$ is compact from $L^{\Psi}(\mathcal{A})$ into $L^{\Psi}(\Sigma)$.
 By the same method of (a) and theorem 2.1 (b), we can find the sequence $\{f_{n}\}_{n\in \mathbb{N}}$ in $L^{\Psi}(\mathcal{A})$, such that $N_{\Psi}(f_{n})=1$, and

$$\int_{\Omega}\Psi(uf_{n})d\mu=\int_{\Omega}\Psi(M_{u}(f_{n}))d\mu
\geq\Psi(\varepsilon_{0})\Psi(\frac{1}{b^2}).$$
Since $\Psi$ is increasing and $|f_{n}-f_{m}|=|f_{n}|+|f_{m}|$, then

$$\int_{\Omega}\Psi(|uf_{n}-uf_{m}|)d\mu=\int_{\Omega}\Psi(|M_{u}(|f_{n}|+|f_{m}|)|)d\mu
\geq\Psi(\varepsilon_{0})\Psi(\frac{1}{b^2}).$$

Thus $\{uf_{n}\}_{\in \mathbb{N}}$ has no convergence subsequence
in $\Psi-$mean convergence. So $\{uf_{n}\}_{\in \mathbb{N}}$ has
no convergence subsequence in norm. This is a contradiction. \ \

c) Let $\varepsilon>0$ and $N_{\varepsilon}=N_{\varepsilon}(\Psi^{-1}(E(\Psi(u)))$. By assumption,
there exist finitely many disjoint $\mathcal{A}-$atoms $\{A_{i}\}^{n}_{i=1}$, such that $N_{\varepsilon}=\cup^{n}_{i=1}A_{i}$.
 Define the operator $T_{\varepsilon}$ on $L^{\Phi}(\Sigma)$, such that $T_{\varepsilon}(f)=E(uf\chi_{N_{\varepsilon}})$, for $f\in L^{\Phi}(\Sigma)$.
  Since  $N_{\varepsilon}\in \mathcal{A}$, $E(uf)$ is $\mathcal{A}-$measurable and  $\mathcal{A}-$measurable functions are constant on  $\mathcal{A}-$atoms,
  we have

$$T_{\varepsilon}(f)=E(uf)\chi_{N_{\varepsilon}}=\sum^{n}_{i=1}E(uf)(A_{i})\chi_{A_{i}}\in L^{\Phi}(N_{\varepsilon}).$$

So $T_{\varepsilon}$ is finite rank.

For $f\in L^{\Phi}(\Sigma)$,
$$T(f)-T_{\varepsilon}(f)=E(uf)-E(uf)\chi_{N_{\varepsilon}}=E(uf)\chi_{\Omega \setminus N_{\varepsilon}}.$$

Thus

$$\int_{\Omega}\Phi(\frac{T(f)-T_{\varepsilon}(f)}{C\varepsilon N_{\Phi}(f)})d\mu=\int_{\Omega}\Phi(\frac{E(uf)\chi_{\Omega \setminus N_{\varepsilon}}}{C\varepsilon N_{\Phi}(f)})d\mu=$$

 $$\leq\int_{\Omega\setminus N_{\varepsilon}}\Phi(\frac{C \Phi^{-1}(E(\Phi(|\frac{f}{N_{\Phi}(f)})|)))\Psi^{-1}(E(\Psi(|u|)))}{C\varepsilon})d\mu$$

$$ \leq\int_{\Omega\setminus N_{\varepsilon}}\Phi(\Phi^{-1}(E(\Phi(|\frac{f}{N_{\Phi}(f)})|)))d\mu\leq\int_{\Omega}E(\Phi(\frac{f}{N_{\Phi}(f)}))d\mu$$

$$=\int_{\Omega}\Phi(\frac{f}{N_{\Phi}(f)})d\mu\leq1.$$

This implies that $N_{\Phi}(T(f)-T_{\varepsilon}(f))\leq C\varepsilon N_{\Phi}(f)$, and so $\|T-T_{\varepsilon}\|<C\varepsilon$. This mean's that, $T$ is limit of a sequence of finite rank operators. So $T$ is compact.

\
\

{\bf Corollary 2.4.} \\(a) If  $(E, \Phi)$ satisfies in
Generalized conditional-type Holder inequality and $\Phi\in
\bigtriangleup'$(globally), then $T$ is compact if and only if
 $N_{\varepsilon}(\Psi^{-1}(E(\Psi(u)))$ consists of finitely many $\mathcal{A}-$atoms, for all
 $\varepsilon>0$.\\

(b) If $\Psi\prec x$(globally) and $(E, \Phi)$ satisfies in
Generalized conditional-type Holder inequality, then then $T$ is
compact if and only if
 $N_{\varepsilon}(\Psi^{-1}(E(\Psi(u)))$ consists of finitely many $\mathcal{A}-$atoms, for all
 $\varepsilon>0$.\\

(c)If $(\Omega, \mathcal{A}, \mu)$ is non-atomic measure
space,$(E, \Phi)$ satisfies in Generalized conditional- type
Holder inequality and $\Phi\in \bigtriangleup'$(globally). Then
$T=EM_{u}$ is a compact operator on $L^{\Phi}(\Sigma)$ if and only
if $T=0$.

{\bf Proof} (b) Since $\Psi\prec x$
 then $\Psi(u)\leq Ku$ for some $K>0$. Suppose that, there exists $\varepsilon_{0}>0$,
such that $N_{\varepsilon}(\Psi^{-1}(E(\Psi(u))))$ doesn't consist
finitely many $\mathcal{A}-$atoms.
 By the same method of (a) and theorem 2.1 (b), we can find the sequence $\{f_{n}\}_{n\in \mathbb{N}}$ in $L^{\Phi}(\mathcal{A})$,
 such that $N_{\Phi}(f_{n})=1$, and

$$\int_{\Omega}\Phi(Kuf_{n})d\mu\geq\int_{\Omega}\Phi(\Psi(u)Kf_{n})d\mu$$$$=\int_{\Omega}E(\Phi(\Psi(u)Kf_{n}))d\mu\geq
\int_{E_n}\Phi(E(\Psi(u))f_{n})d\mu $$$$\geq\varepsilon_{0}K,$$
where $K\varepsilon_0>1$. Since $\Phi$ is increasing and
$|f_{n}-f_{m}|=|f_{n}|+|f_{m}|$, then

$$\int_{\Omega}\Phi(|Kuf_{n}-Kuf_{m}|)d\mu=\int_{\Omega}\Phi(|KM_{u}(|f_{n}|+|f_{m}|)|)d\mu
\geq K\varepsilon_{0}.$$

Thus $\{KM_u(f_{n})\}_{\in \mathbb{N}}$ has no convergence
subsequence in $\Phi-$mean convergence. So $\{M_u(f_{n})\}_{\in
\mathbb{N}}$ has no convergence subsequence in norm. This is a
contradiction. \

 \

\
 In the next theorem we use the method that is used in \cite{her}.\\

{\bf Theorem 2.5.} If $\mathcal{A}\neq\Sigma$, then $\sigma(T)=ess
range(E(u))\cup\{0\}$.

\
\

{\bf Proof.} Since $\mathcal{A}\neq\Sigma$,  $L^{\Phi}(\mathcal{A})\varsubsetneq L^{\Phi}(\Sigma)$. Hence $T=EM_{u}$ isn't surjective and so
$0\in\sigma(T)$.
Let $\lambda\notin ess range(E(u))$, $\lambda\neq0$. We show that $T-\lambda I$ is invertible. If $Tf-\lambda f=0$, then $E(uf)=\lambda f$. So $f$ is $\mathcal{A}-$measurable. Thus
$(E(u)-\lambda)f=E(uf)-\lambda f=0$. Since $\lambda\notin ess range(E(u))$, then $E(u)-\lambda\geq\varepsilon$ a.e for some $\varepsilon>0$. So
$f=0$ a.e. This implies that $T-\lambda I$ is injective.
Now we show that $T-\lambda I$ is surjective. Let $g\in L^{\Phi}(\Sigma)$. We can write

$$g=g-E(g)+E(g), \ \ \ g_{1}=g-E(g),\ \ \ g_{2}=E(g).$$

Since $N_{\Phi}(E(g))\leq N_{\Phi}(g)$, then $g_{2}\in L^{\Phi}(\mathcal{A})$ and  $g_{1}\in L^{\Phi}(\Sigma)$, $(E(g_{1})=0$.
Let $$f_{1}=\frac{\lambda g_{1}+T(g_{2})}{\lambda(E(u)-\lambda)}, \ \ \ f_{2}=\frac{-g_{2}}{\lambda}.$$
Since $\lambda\notin ess range(E(u))$, then $E(u)-\lambda\geq\varepsilon$ a.e for some $\varepsilon>0$. So $\|\frac{1}{E(u)-\lambda}\|_{\infty}\leq\frac{1}{\varepsilon}$. Thus $f_{2}\in L^{\Phi}(\mathcal{A})$,  $f_{1}\in L^{\Phi}(\Sigma)$ and $f=f_{1}+f_{2}\in  L^{\Phi}(\Sigma)$. Direct computation shows that $T(f)-\lambda f=g$. This implies that $T-\lambda I$ is invertible and so $\lambda \notin \sigma(T)$.

\
\

Conversely, let  $\lambda\notin \sigma(T)$. Define linear transformation $S$ on $L^{\Phi}(\Sigma)$ as follows

$$Sf=\frac{Tf-f(E(u)-\lambda)}{\lambda(E(u)-\lambda}, \ \ \ f\in L^{\Phi}(\Sigma).$$

If  $\lambda\notin ess range(E(u))$, then $\|\frac{1}{E(u)-\lambda}\|_{\infty}\leq\frac{1}{\varepsilon}$ for some $\varepsilon>0$. So
$$N_{\Phi}(Sf)\leq N_{\Phi}(\frac{Tf}{\lambda(E(u)-\lambda})+N_{\Phi}(\frac{f}{\lambda})$$

$$\leq(\frac{\|T\|}{\lambda \varepsilon}+\frac{1}{\varepsilon})N_{\Phi}(f).$$
Thus $S$ is bounded $L^{\Phi}(\Sigma)$.
If $S$ is bounded on $L^{\Phi}(\Sigma)$, then for $f\in L^{\Phi}(\mathcal{A})$
$Sf=\alpha f=M_{\alpha}f$, where $\alpha=\frac{1}{Eu-\lambda}$. Thus multiplication operator $M_{\alpha}$ is bounded on $L^{\Phi}(\mathcal{A})$.
This implies that $\alpha\in L^{\infty}(\mathcal{A}$ and so there exist some $\varepsilon>0$ such that $E(u)-\lambda=\frac{1}{\alpha}\geq\varepsilon$ a.e.
 This mean's that $\lambda\notin ess range(E(u))$. Also, we have

$$S\circ (T-\lambda I)=(T-\lambda I)\circ S=I.$$
Thus $(T-\lambda I)^{-1}=S$ and so  $\sigma(T)=ess range(E(u))\cup\{0\}$.
\\

Let $\mathfrak{B}$ be a Banach space and $\mathcal{K}$ be the set
of all compact operators on $\mathfrak{B}$. For $T\in
L(\mathfrak{B})$, the Banach algebra of all bounded linear
operators on $\mathfrak{B}$ into itself, the essential norm of $T$
means the distance from $T$ to $\mathcal{K}$ in the operator norm,
namely $\|T\|_e =\inf\{\|T - S\| : S \in\mathcal{K}\}$. Clearly,
$T$ is compact if and only if $\|T\|_e= 0$. Let $X$ and $Y$ be
reflexive Banach spaces and $T\in L(X,Y)$. It is easy to see that
$\|T\|_{e}=\|T^*\|_{e}$. As is seen in \cite{sha}, the essential
norm plays an interesting role in the compact problem of concrete
operators.\\

 In the sequel we present an upper bound for
essential norm of $EM_u$ on Orlicz space $L^{\Phi}(\Sigma)$.
\\
\vspace*{0.3cm} {\bf Theorem 2.6.} Let $EM_{u}: L^{\Phi}(\Omega,
\Sigma, \mu)\rightarrow L^{\Phi}(\Omega, \Sigma, \mu)$ is bounded
and $(E, \Phi)$ satisfies in Generalized conditional-type Holder
inequality. Let $\beta=\inf\{\varepsilon>0:N_{\varepsilon}$
consists of finitely many $\mathcal{A}$-atoms$\}$, where
$N_{\varepsilon}=N_{\varepsilon}(\Psi^{-1}(E(\Psi(u))))$. Then

\

$$\|EM_{u}\|_{e}\leq \beta.$$

{\bf Proof} Let $\varepsilon>0$. Then $N_{\varepsilon+\beta}$
consist of finitely many $\mathcal{A}$-atoms. Put
$u_{\varepsilon+\beta}=u\chi_{N_{\varepsilon+\beta}}$ and
$EM_{u_{\varepsilon+\beta}}$. So $EM_{u_{\varepsilon+\beta}}$ is
finite rank and so compact. By the same method that is used in
theorem 2.3(c) we have

$$\|EM_u\|_{e}\leq \|EM_u-EM_{u_{\varepsilon+\beta}}\|\leq \beta+\varepsilon.$$

This implies that $\|EM_{u}\|_{e}\leq \beta$.

\end{document}